\documentclass[final]{siamltex}
\usepackage{amsmath}
\usepackage{amsfonts}
\usepackage{graphicx}

\def\bJ{\mathbf J}

\def\b0{\mathbf 0}

\title{Phase-space analysis of large ODE systems using a
low-dimensional conservation law
}

\author{Neil~V.~Budko\thanks{Numerical Analysis, DIAM,
        Delft University of Technology, 2628 CD Delft, Mekelweg 4,
        The Netherlands ({\tt n.v.budko@tudelft.nl})}
	\and Fred~J.~Vermolen\footnotemark[1]}

\begin{document}

\maketitle

\begin{abstract}
Simultaneous deterministic and weakly stochastic dynamics of multiple populations 
described by a large system of ODE's is considered in the phase space of population 
sizes and ODE's parameters. We show that many practically interesting problems
can be formulated as a low-dimensional phase-space conservation law and
solved either explicitly or with simple iterative methods. 
In particular, we consider: non-interacting populations with unbounded and 
logistic growth, populations with randomized and biased migration, populations competing 
for a resource, coexisting species, and populations with phase-space interactions. 
The method provides an alternative to Monte Carlo simulations and 
may be useful in the fast analysis of biological data and/or 
removal of deterministic trends.
\end{abstract}

\pagestyle{myheadings}
\thispagestyle{plain}
\markboth{N.~V.~BUDKO AND F.~J.~VERMOLEN}{PHASE-SAPCE ANALYSIS IN POPULATION DYNAMICS}

\section{Introduction}
Many natural phenomena
can be described as a large number of simultaneously developing 
`populations'. The various coexisting species is the most obvious case.
However, even a single plant consists of many cells that in turn contain 
significant numbers of mitochondria -- small organelles capable of 
multiplying on their own. Thus, considering the mitochondria of 
a single cell to be the elementary population, one can view 
each individual plant as a population of 
populations, or a {\it multi-population}. 
On another level, viewing a single plant as an elementary population of cells, a
field of plants is a population of populations, etc. (see, e.g., Merchant et~al.~1960; 
Carrie et~al.~2012).

Although, classical population dynamics is a very well developed subject
(see, e.g., Shonwinkler and Herod~2009; Begon et~al.~2006; Newman et~al. 2014), its main concern
has traditionally been the growth of a single or at most of a few populations simultaneously. 
The question of dynamics of dispersing populations has 
led to the development of the metapopulation theory (e.g., Levin~1969; Eriksson et~al.~2014)
that treats populations with spatially disjoint habitats as a single 
metapopulation, and reduces the mathematical description to a single ODE. 
Alternatively, it is sometimes advantageous to consider the space-time dynamics of all members 
of all populations simultaneously, perhaps, subdividing them in a few species as well. 
This leads to spatially-resolved PDE models, such as the Fischer 
equation (Fischer~1930), that have produced many valuable insights over 
the years (see, e.g., Edelstein-Keshet~2005).

Mathematically, all problems under consideration in this paper represent large
systems of coupled ODE's. 
We resort to the phase-space analysis for the following reasons. 
Firstly, studying multiple populations with large systems of ODE's 
or spatially-resolved PDE's does not always 
provide explicit answers to the pertaining questions. 
For example, in agriculture and ecology (e.g., Hautala and Hakoj{\"a}rvi~2011; 
Lin et~al.~2013; Xu et~al.~2014; Velazquez et~al.~2014), where the growth of
an individual plant or species is typically described by a logistic or a reaction-type ODE,
the actual question of interest or the measured data is often the histogram/distribution 
of plant or population sizes. Such a distribution is then typically obtained with a
Monte-Carlo simulation, where the individual ODE's are solved
over a set of points in the range of their parameters. 
It is a well-known fact that the choice 
of these points (sampling) may significantly influence the quality of the 
prediction of the corresponding distribution function.
Secondly, although the phase-space formulation
results in a PDE rather than ODE, we show that many practically interesting
systems of ODE's feature special types of coupling
corresponding to low-dimensional phase-space problems, 
where the solution is simple to obtain. 

To avoid further confusion we mention
that the term `dimension' refers here to the dimension of the underlying 
manifold (number of phase-space variables) not the dimension of the 
phase space itself. Formally, 
the phase space is infinite-dimensional in our formulation, 
since the state of the system at any given time is represented by a continuous
function of phase-space variables, not just one point as with the usual dynamical systems 
described by a few ODE's, such as the predator-prey system (see, e.g., Edelstein-Keshet~2005). 

To predict the time evolution of the distribution function 
we employ the phase-space conservation law, where the 
mathematical form of the phase-space current is determined by the dynamic 
equations (coupled or uncoupled ODE's) of individual populations. This approach 
is a systematic extension of the Lifshitz-Slyozov-Wagner method from materials science
(Lifshitz and Slyozov~1961; Wagner~1961; Kampman and Wagner~1991;
Myhr and Grong~2000; Collet and Goudon~2000; Laurencot~2002).
Conceptually it may be viewed as a version of the
Fokker-Planck equation tuned, in our case, to the specific needs of the multi-population 
dynamics. The present approach is also similar in spirit to the 
Liouville equation of Hamiltonian mechanics, except that the population dynamics
is not necessarily Hamiltonian. The closest analogue in the field of population 
dynamics is, probably, the state-space method discussed in (Newman et~al. 2014).

In the next section we present the general formulation of the problem and
elucidate the meaning and the intended purpose of our multi-population 
distribution functions.
Next, we illustrate the benefits and limitations of the proposed framework
on seven population models relevant to plant biology, ecology, and migration studies. 
We focus on problems that have either completely explicit solutions or
can be solved with a simple iterative algorithm.
In particular, we consider populations with unlimited (exponential) 
growth, populations whose growth is limited by the carrying capacity of the medium 
(logistic equation), populations/species competing for resources 
(e.g. oxygen, food, etc.), populations with migrating members, 
and populations whose rate of growth depends on the 
distribution function. The mathematical analysis of the resulting PDE's 
is mostly available in the literature, apart from the case of the coexisting species,
which is studied in some details in the corresponding section.

\section{General formulation}
Let the growth of $P$ populations be governed by the following system of 
ODE's:
\begin{align}
\label{eq:GenODE}
\begin{split}
&\frac{d n_{i}}{dt} = f_{i}(n_{1},\dots,n_{P},\alpha_{1},\dots,\alpha_{k},t);
\;\;\;i=1,\dots,P,
\\
&\frac{d \alpha_{j}}{dt} = g_{j}(t),\;\;\;j=1,\dots,k,
\end{split}
\end{align}
where $f_{i}$ and $g_{j}$ are differentiable functions of the indicated variables,
and the variables $\alpha_{j}$, $j=1,\dots,k$, represent the various dynamical parameters
of the ODE's, such as the rates of birth and death. 

In general, the phase-space 
formulation of this problem is even more difficult than the original system of ODE's.
For example, if all rate functions $f_{i}$ were, indeed, different, then the associated 
phase-space manifold is $(P+k)$-dimensional, i.e., the size of each population 
requires a separate coordinate. Hence, from the practical point of view, 
the phase-space analysis of large ODE systems
only makes sense if the number of dimensions of the phase-space manifold is much smaller than 
$P$. The possibility of such a low-dimensional phase-space formulation depends on
the type of the rate functions $f_{i}$. Here we shall limit ourselves to 
a specific $(k+1)$-dimensional manifold, which in its turn limits the class of 
admissible dynamic systems.

Our goal is the distribution function 
$u(n,\alpha_{1},\dots,\alpha_{k},t): {\mathbb R}^{k+2}\rightarrow{\mathbb R}$ with 
the following properties:
\begin{align}
\label{eq:GenUProp1}
&\int\dots\int\int u(n,\alpha_{1},\dots,\alpha_{k},t)\,dn\,d\alpha_{1}\dots d\alpha_{k} = P,
\\
\label{eq:GenUProp2}
&\int\dots\int\int n u(n,\alpha_{1},\dots,\alpha_{k},t)\,dn\,d\alpha_{1}\dots d\alpha_{k} = 
\sum_{i=1}^{P}n_{i}(t)=N(t),
\\
\label{eq:GenUProp3}
&u(n,\alpha_{1},\dots,\alpha_{k},0)=u_{0}(n,\alpha_{1},\dots,\alpha_{k}),
\end{align}
where $u_{0}$ is the given initial distribution of the populations over the
$(k+1)$-dimensional phase space $\Omega$.

The conservation of the total number of populations
means that the time-evolution of the distribution function $u$ is governed by the 
generalized continuity equation:
\begin{align}
\label{eq:GenContinuity}
\frac{\partial u}{\partial t} 
+ \frac{\partial J_{n}}{\partial n} +  
\sum_{j=1}^{k}\frac{\partial J_{\alpha_{j}}}{\partial \alpha_{j}} = 0,
\end{align}
where the phase-space current density $\bJ$ has the form
\begin{align}
\label{eq:GenJ}
\bJ=\langle J_{n},J_{\alpha_{1}},\dots, J_{\alpha_{k}} \rangle
=\langle uv_{n},uv_{\alpha_{1}},\dots, uv_{\alpha_{k}} \rangle.
\end{align}
Since, e.g., $v_{n}=dn/dt$, by definition, the phase-space `velocities' 
are to be deduced from the corresponding ODE's:
\begin{align}
\label{eq:GenVelocity}
\begin{split}
&v_{n} = \tilde{f}(n,\alpha_{1},\dots,\alpha_{k},t),
\\
&v_{\alpha_{j}} = g_{j}(t),\;\;\;j=1,\dots,k.
\end{split}
\end{align}
where $\tilde{f}$ is the phase-space representation of the functions $f_{i}$.
The success of this $(k+1)$-dimensional approach depends on the existence of 
a one-to-one correspondence between $f_{i}$ and 
$\tilde{f}$, i.e., one should be able to express $f_{i}$ purely in terms 
of the chosen set of phase-space coordinates $\{n,\alpha_{1},\dots,\alpha_{k}\}$ and 
time $t$. Below we describe several types of such $f_{i}$'s. 

The first type
corresponds to the uncoupled system of ODE's and a uniform (independent of index $i$) 
rate function:
\begin{align}
\label{eq:GenUncoupled}
f_{i}(n_{1},\dots,n_{P},\alpha_{1},\dots,\alpha_{k},t)
=f(n_{i},\alpha_{1},\dots,\alpha_{k},t).
\end{align}
For example, the polynomial rate functions
\begin{align}
\label{eq:RatesUncoupled}
f_{i}=\alpha_{i}+\beta_{i}n_{i}+\gamma_{i}n_{i}^{2}+\dots,
\end{align}
correspond to the phase-space velocity
\begin{align}
\label{eq:VelocitiesUncoupled}
v_{n}=\tilde{f}(n,\alpha,\beta,\gamma,\dots,t)=\alpha+\beta n+\gamma n^{2}+\dots.
\end{align}
Such uncoupled systems result in linear low-dimensional phase-space Cauchy 
problems that can be solved completely explicitly.

The second type of $f_{i}$ is also uniform and features a uniform coupling of the form:
\begin{align}
\label{eq:GenCoupled}
f_{i}(n_{1},\dots,n_{P},\alpha_{1},\dots,\alpha_{k},t)
=f(n_{i},N(t),\alpha_{1},\dots,\alpha_{k},t),
\end{align}
where $N(t)$ is given by (\ref{eq:GenUProp2}). For example, if the coupling has the 
form
\begin{align}
 \label{eq:UniformCoupling}
f_{i}=\alpha_{i}\sum_{j}n_{j}=\alpha_{i} N(t),
\end{align}
then the corresponding phase-space velocity is
\begin{align}
 \label{eq:UniformCouplingPhaseSpace}
v_{n}=\tilde{f}(\alpha,t)=\alpha \iint n u(n,\alpha,t)\,dn\,d\alpha.
\end{align}
The resulting weakly nonlinear low-dimensional phase-space
problems can be solved with simple iterative methods. 

However, the important general case, where the coupling has the form:
\begin{align}
 \label{eq:NonUniformCoupling}
f_{i}=\sum_{j}\alpha_{ij}n_{j},
\end{align}
would require a separate phase-space coordinate 
for each population and is, therefore, impractical.

A further generalization of (\ref{eq:GenCoupled}) that admits an explicit phase-space 
representation has the form
\begin{align}
\label{eq:PhaseSpaceCoupled}
f_{i}(n_{1},\dots,n_{P},\alpha_{1},\dots,\alpha_{k},t)
=f(n_{i},F[u],\alpha_{1},\dots,\alpha_{k},t),
\end{align}
where $F[u]$ is a functional of the distribution function $u$. Below
we consider one such problem that leads to the Burgers equation in 
the phase space.

We conclude this section with a few words about the interpretation of the
proposed phase-space picture. Since the population dynamics considered here is mainly deterministic, 
the functions $u_{0}$ and $u$ are not probability density functions, 
although, upon suitable normalization, they can be 
interpreted as such for a randomly selected population.

The relation of the distribution function $u$ to the populations whose dynamics is
governed by (\ref{eq:GenODE}) is similar to the relation between the data, their histogram,
and the corresponding probability density function. It is even more similar to the 
continuous description of discrete natural phenomena common in the classical 
macroscopic theories of mechanics, fluid motion, and electromagnetism.

It is well-known that one has to be careful when choosing the appropriate bin size for 
a histogram. Take it too small and the histogram breaks down into a meaningless 
collection of separated bars each containing only one data point, i.e., 
each having height equal to one.
Similarly, in macroscopic continuum physics the basic assumption is the
so-called elementary volume, which while being relatively small contains a sufficient number 
of particles, so that the 
concept of number density makes sense. Generally, the continuum hypothesis and the
corresponding equations are not considered to be valid at scales smaller than the 
elementary volume.

The intended purpose and the desired property of the distribution function $u$
is that it approximates the number $\delta P$ of
populations with their parameters within a subregion $\delta \Omega$ of the 
phase space, i.e.,
\begin{align}
\label{eq:GenInterpretation}
\delta P(\delta\Omega,t)\approx
\int \dots \int_{\delta \Omega} u(n,\alpha_{1},\dots,\alpha_{k})\,d\Omega.
\end{align}
For this interpretation to be valid (i.e., reasonably precise) the region $\delta \Omega$ 
must be sufficiently large. Obviously, the interpretation is exact, 
if $\delta \Omega$ coincides with the whole of $\Omega$. The acceptable lower bound on 
$\delta\Omega$, in general, will depend on both the number of populations $P$ and
the smoothness/extent of the distribution $u$.

In what follows we apply our approach to a series of
models of increasing complexity. We restrict ourselves to the so-called 
autonomous case, where the dynamic parameters $\alpha_{j}$ are constant in time.
On the one hand, extension to the non-autonomous case of time-dependent $\alpha_{j}$'s
is trivial. On the other hand, the practical aspects of problems requiring such an extension 
deserve a detailed investigation, which is beyond the scope of the present paper. 
an $n$-dimensional 

\section{Examples of problems with low-dimensional conservation laws}
\subsection{Populations with unlimited growth}
Let a very large number of populations be growing in accordance with the following equations:
\begin{align}
\label{eq:nUnlimited}
\frac{dn_{i}}{dt}=\alpha_{i} n_{i},\;\;\;n_{i}(0)> 0,
\;\;\;i=1,2,\dots P,
\end{align}
where $\alpha_{i}> 0$ are the growth constants. The variable $n_{i}$ 
will be referred to as the size of the $i$-th population, but, in practice, it may denote 
many different things. For example, it can be the length, volume, or the biomass of a plant. 
It can also be a continuous approximation of the number of animals within a 
habitat, or the number of cells in a plant, or the number of mitochondria in a cell. 
The units of time are also completely arbitrary here.

As was explained in the previous section, instead of solving all these ODE's we consider the 
time evolution of the distribution function $u(n,\alpha,t)$, whose
integral over the size variable $n$ and the growth constant variable
$\alpha$ represents the total number of populations, which is conserved:
\begin{align}
\label{eq:uDoubleIntegral}
\int\limits_{0}^{\infty}\int\limits_{0}^{\infty}
u(n,\alpha,t)\,dn\,d\alpha = P.
\end{align}
Integrating the distribution over $\alpha$ only,
\begin{align}
\label{eq:uAlphaIntegral}
\int\limits_{0}^{\infty}
u(n,\alpha,t)\,d\alpha = \rho(n,t),
\end{align}
one obtains a continuous approximation of the distribution of 
populations over their sizes at time $t$. 

We assume that the sufficiently smooth continuous approximation of the initial
distribution $u(n,\alpha,0)=u_{0}(n,\alpha)$ is given.
The latter is interpreted as a two-dimensional distribution of populations 
over their initial sizes $n_{0}$ and their growth constants $\alpha$.

In the present case the phase-space current has at most two components
\begin{align}
\label{eq:J}
\bJ=\langle uv_{n},uv_{\alpha} \rangle,
\end{align}
and the corresponding `velocities' are easily deduced from the dynamic 
equation (\ref{eq:nUnlimited}) as: 
$v_{n}=\alpha n$ and $v_{\alpha}=0$.
Since $\bJ$ has only one nonzero component, 
the continuity equation (\ref{eq:GenContinuity}) reduces to the following Cauchy
problem:
\begin{align}
\label{eq:uUnlimited}
\frac{\partial u}{\partial t} + \alpha n\frac{\partial u}{\partial n} + \alpha u = 0,
\;\;\;u(n,\alpha,0)=u_{0}(n,\alpha).
\end{align}
The solution of this problem, obtained with the method of characteristics, is
\begin{align}
\label{eq:uSolution}
u(n,\alpha,t)=e^{-\alpha t}u_{0}(ne^{-\alpha t},\alpha).
\end{align}
\begin{figure}[t]
\begin{center}
\hspace*{-0.5cm}
\includegraphics[width=13cm]{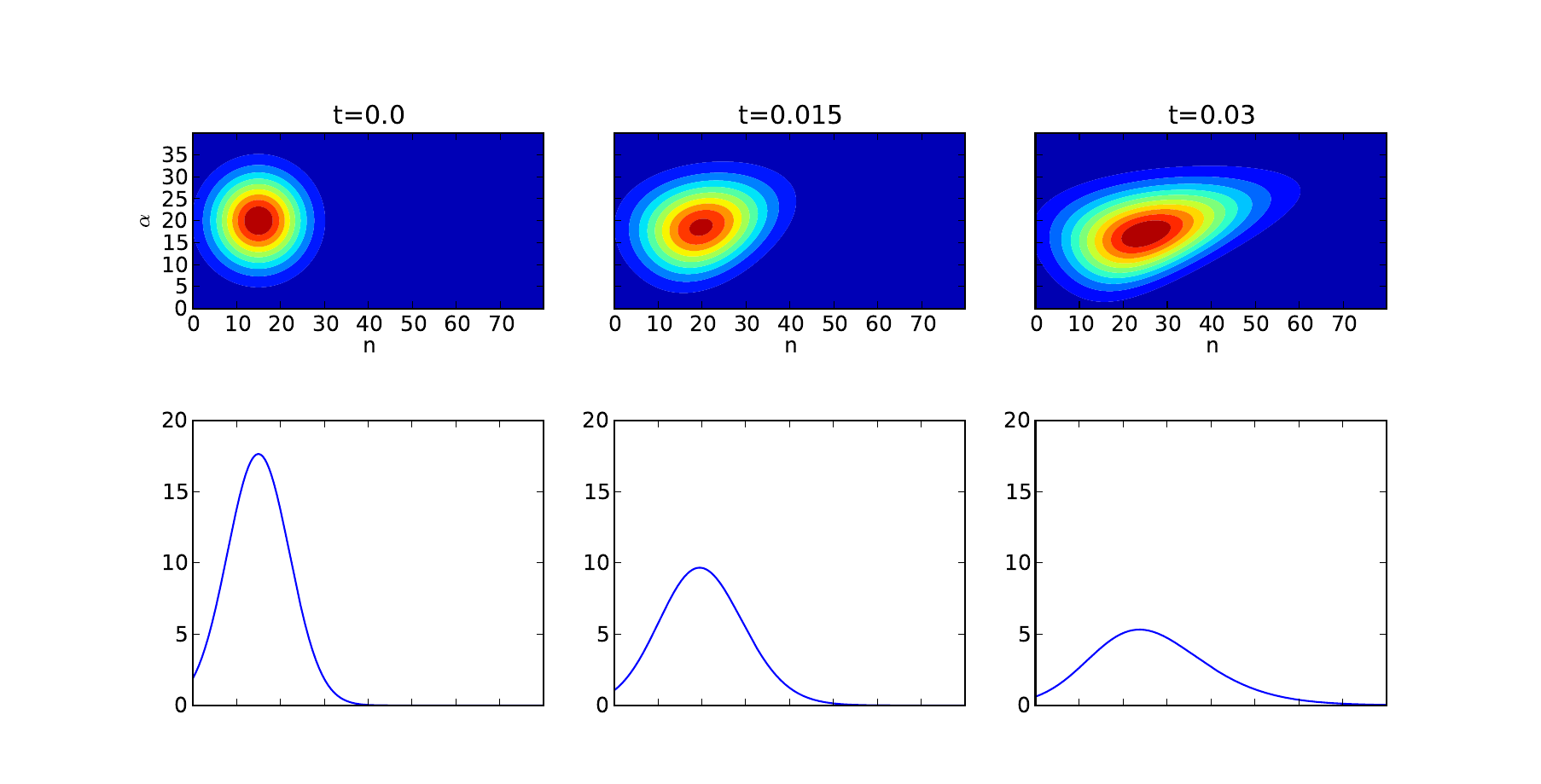}
\end{center}
\caption{
\label{fig:Fig1}
Exponentially growing populations.  
Top row -- snapshots of $u(n,\alpha,t)$ at $t=0.0, 0.015$, and $0.03$; Bottom 
row -- corresponding numerically integrated $\rho(n,t)$ showing the distribution of
populations over their sizes (the horizontal, $n$-axis, of these plots 
is shared with the images above).
}
\end{figure}
An example of the time evolution of $u(n,\alpha,t)$ is shown in Fig.~\ref{fig:Fig1},
where the initial distribution $u_{0}(n,\alpha)$ is a Gaussian centered at $n_{\rm c}=15$ and 
$\alpha_{\rm c}=20$, which was set to zero for $\alpha,n < 0$ (see the top-left image). 
The size distributions of populations $\rho(n,t)$ 
(bottom row of Fig.~\ref{fig:Fig1}) were obtained by numerical integration
over $\alpha$ and demonstrate a progressively fattening tail due to the rapidly
growing populations corresponding to the upper part of $u(n,\alpha,t)$,
i.e., larger growth rate constants $\alpha$. An interesting feature in 
the evolution of $u(n,\alpha,t)$ that becomes
more pronounced at larger $t$'s, 
is the gradual downwards shift of the maximum
of $u(n,\alpha,t)$ with time. This is due to the fact that the populations
with smaller $\alpha$'s tend to grow more evenly (synchronously) even if their 
initial $n$'s are different, whereas, populations with larger $\alpha$'s 
and slightly different initial $n$'s rapidly grow apart.

\subsection{Populations with logistic growth}
Let the dynamics of populations be governed by the logistic equations
\begin{align}
\label{eq:Logistic}
\frac{dn_{i}}{dt}=\gamma_{i} n_{i}(k_{i}-n_{i}),\;\;\;n_{i}(0)> 0,
\;\;\;i=1,2,\dots,
\end{align}
where $\gamma_{i} \geq 0$ is a rate constant, and $k_{i}>0$ is a constant determining
the maximum sustainable population. The distribution function 
$u(n,\gamma,k,t)$ over the three-dimensional $(n,\gamma,k)$ phase space satisfies
the following Cauchy problem:
\begin{align}
\label{eq:uLogistic}
\frac{\partial u}{\partial t} + \gamma n(k-n)\frac{\partial u}{\partial n}
+ \gamma(k-2n) u = 0,\;\;\;u(n,\gamma,k,0)=u_{0}(n,\gamma,k).
\end{align}
Here too the solution can be obtained with the method of characteristics. 
Since it is more involved than in the previous case, we shall provide the 
intermediate steps. The characteristic equations of (\ref{eq:uLogistic}) are
\begin{align}
\label{eq:CharacteristicLogistic}
\frac{dt}{ds}=1,\;\;\;
\frac{dn}{ds}=\gamma n(k-n),\;\;\;
\frac{du}{ds}=-\gamma (k-2n)u.
\end{align}
This leads to the following relations:
\begin{align}
\label{eq:CharacteristicLogisticSolutions}
C_{0}=\frac{n}{\vert k-n\vert}e^{-\gamma kt},
\;\;\;
u=\frac{C_{1}}{n\vert k-n\vert}.
\end{align}
We shall restrict ourselves to the segment $0\leq n<k$, so that $\vert k-n\vert=k-n$. 
The general form of the solution to (\ref{eq:uLogistic}) in this case is
\begin{align}
\label{eq:uLogisticGeneral}
u=\frac{f(p)}{n(k-n)},\;\;\;p=\frac{n}{k-n}e^{-\gamma kt}.
\end{align}
The function $f(p)$ should be chosen in such a way that the initial condition
is satisfied. To this end we introduce two smooth functions of $p$ that at 
$t=0$ behave as follows:
\begin{align}
\label{eq:pAt0}
\left.\frac{kp}{p+1}\right\vert_{t=0}=n,\;\;\;\;\;
\left.\frac{k^{2}p}{(p+1)^{2}}\right\vert_{t=0}=n(k-n).
\end{align}
The initial condition from (\ref{eq:uLogistic}) is satisfied if these functions are used as
\begin{align}
\label{eq:uLogisticInitial}
\left.\frac{k^{2}p}{n(k-n)(p+1)^{2}}\;\;
u_{0}\left(\frac{kp}{p+1},\gamma,k\right)\right\vert_{t=0}=u_{0}(n,\gamma,k).
\end{align}
Thus, the solution of the Cauchy problem (\ref{eq:uLogistic}) on the strip 
$0\leq n<k$ is
\begin{align}
\label{eq:uLogisticSolution0<N<K}
u(n,\gamma,k,t)=\frac{k^{2}e^{-\gamma kt}}{\left(ne^{-\gamma kt}+k-n\right)^{2}}
\;\;u_{0}\left(\frac{kne^{-\gamma kt}}{ne^{-\gamma kt}+k-n},\gamma,k\right),
\end{align}
and for $n=k$, we obviously have $u(k,\gamma,k,t)=u_{0}(k,\gamma,k)$.

\begin{figure}[t]
\begin{center}
\hspace*{-0.5cm}
\includegraphics[width=13cm]{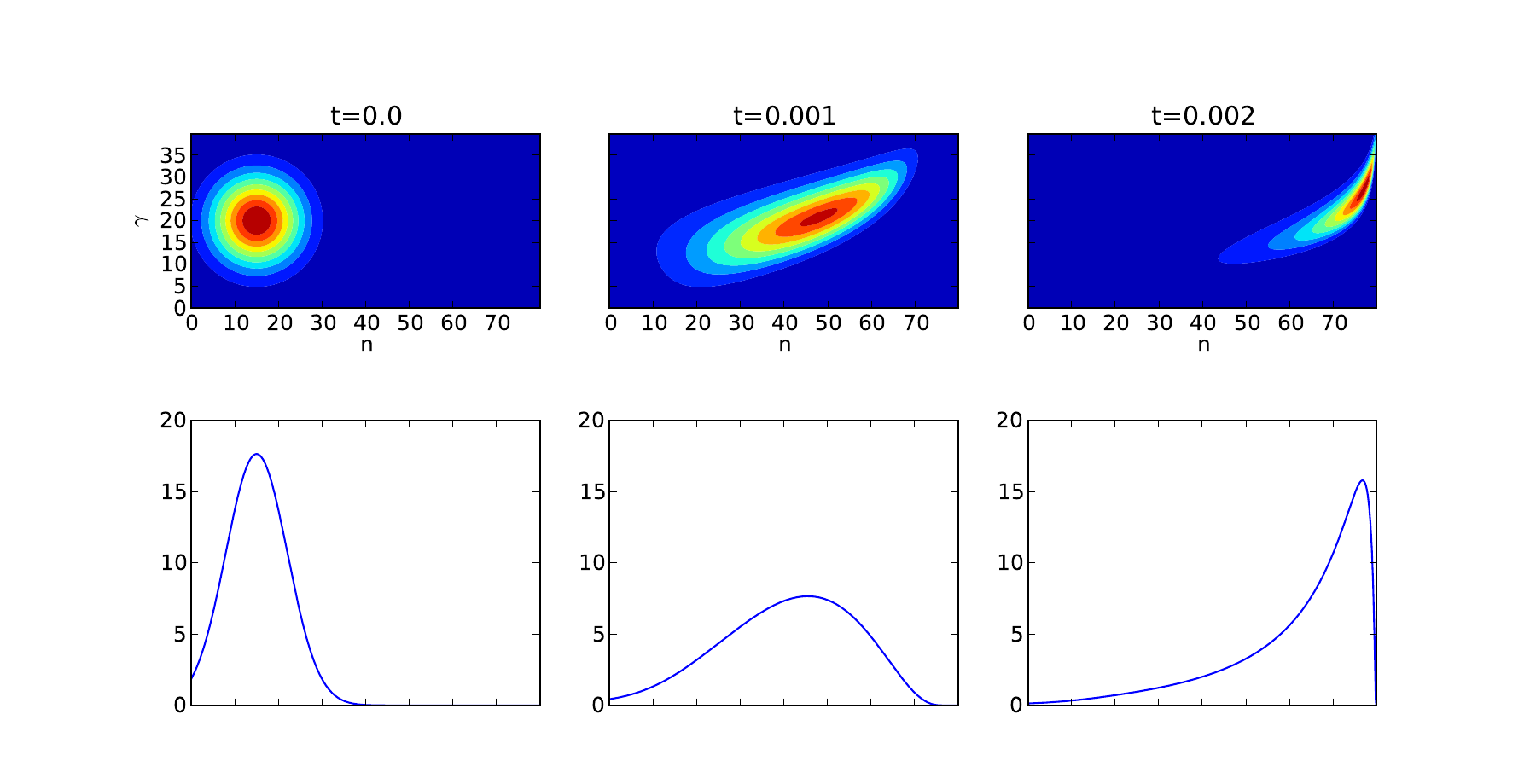}
\end{center}
\caption{
\label{fig:Fig2}
Logistic growth with fixed $k=80$. 
Top row -- snapshots of $u(n,\gamma,k,t)$ at $t=0.0, 0.001$, and $0.002$; Bottom 
row -- corresponding numerically integrated $\rho(n,t)$ showing the distribution of
populations over their sizes.
}
\end{figure}
An example of the time evolution of $u(n,\gamma,k,t)$ for a fixed $k=80$ is shown in 
Fig.~\ref{fig:Fig2},
where the initial distribution $u_{0}(n,\gamma,80)$ (top left) is the same Gaussian 
as in the example of Fig.~\ref{fig:Fig1}. 
As $t$ grows the corresponding size distributions $\rho(n,t)$ (bottom row) 
tend to a distribution localized at the right boundary $n=k$.
This time the maximum of the distribution $u(n,\gamma,k,t)$ shifts upwards (eventually),
meaning that for larger $t$'s the populations with larger growth rates $\gamma$
will dominate, as they will rapidly accumulate around the maximum attainable size.

\subsection{Populations competing for a finite resource}
Consider populations competing for a finite resource $c$ that controls their growth rate.
As an example one can think of a batch of seeds germinating in a sealed container.
Seeds need oxygen to germinate and there is only a finite amount of it in a 
sealed container (Bewley et~al.~2013, van~Duijn and Koenig~2001, van~Asbrouck and Taridno~2009,
Budko et al.~2013). Hence, there will be some sort of passive competition for this 
oxygen among the seeds, but none of the seeds can actually `win', as
all oxygen will be depleted at the end anyway, and all the seeds will suffocate 
(stop their growth).
Although this fate is inevitable, the questions about the time evolution of the distribution 
of seed sizes and the corresponding oxygen consumption curve have some practical 
interest to them, since the initial stages in the growth of the seed are indicative of
the plant vitality.

We further assume that the control substance is absorbed equally well 
by all members of all populations
(with the same reaction rate coefficient $\gamma$).
However, the overall rate of consumption of $c$ is proportional to the 
concentration of $c$ (or partial pressure, if the control substance is a gas).
Populations convert $c$ into energy and biomass and grow. 
The growth of a population slows down and eventually stops if the supply of $c$ drops to 
some minimal value $c_{\rm min}$, which is here taken to be zero for simplicity. 
The following coupled system of reaction-type equations describes this situation:
\begin{align}
\label{eq:CN}
\begin{split}
\frac{dc}{dt}&=-\gamma c \sum_{i}n_{i},\;\;\;c(0)=c_{0},\gamma>0,
\\
\frac{dn_{i}}{dt}&=\beta_{i} n_{i}c,\;\;\;n_{i}(0),\beta_{i}> 0,\;\;\;i=1,2,\dots
\end{split}
\end{align}
Integrating the first equation and substituting the result in the second equation
we arrive at:
\begin{align}
\label{eq:nCompeting}
\frac{dn_{i}}{dt}&=\beta_{i} n_{i}c_{0}e^{-\gamma \tilde{N}(t)},
\;\;\;i=1,2,\dots
\end{align}
where 
\begin{align}
\label{eq:tildeN}
\tilde{N}(t) = \int\limits_{0}^{t}\sum_{i}n_{i}\,dt'
\end{align}
The integrand above is the total size
$N(t)$ of all populations that can be expressed via the phase-space distribution as
\begin{align}
\label{eq:N}
N(t)=\int\limits_{0}^{\infty}\int\limits_{0}^{\infty}
nu(n,\beta,t)\,dn\,d\beta.
\end{align}
Hence, eq.~(\ref{eq:tildeN}) can be rewritten 
in terms of the distribution as well
\begin{align}
\label{eq:TotalPop}
\tilde{N}(t)=\int\limits_{0}^{t}N(t')\,dt'=\int\limits_{0}^{t}\int\limits_{0}^{\infty}\int\limits_{0}^{\infty}
nu(n,\beta,t')\,dn\,d\beta\,dt'
\end{align}
Thus, the phase-space dynamics of the distribution function will be governed by
the weakly nonlinear equation:
\begin{align}
\label{eq:uCompeting}
\frac{\partial u}{\partial t}
+c_{0}\beta ne^{-\gamma \tilde{N}(t)}
\frac{\partial u}{\partial n} + c_{0}\beta u e^{-\gamma \tilde{N}(t)} = 0, \;\;\;u(n,\beta,0)=u_{0}(n,\beta),
\end{align}
The questions of existence and uniqueness of solutions to a 
slightly more general problem of this type 
have been discussed by Collet and Goudon~(2000).
Mathematically, the present problem is not as challenging, since for all $n$ and $\beta$
the velocity of the phase-space flow is positive and goes to zero
exponentially in time, meaning that $\partial u/\partial t$ approaches zero as
$t\rightarrow\infty$.

The method of characteristics yields the following implicit solution:
\begin{align}
\label{eq:uCompetingSolution}
u(n,\beta,t)=e^{-\beta \xi(t)}u_{0}(ne^{-\beta\xi(t)},\beta),
\end{align}
where
\begin{align}
\label{eq:Xi}
\xi(t)&=\int_{0}^{t}c(\tau)\,d\tau,
\\
c(\tau)&= c_{0}e^{-\gamma\tilde{N}(\tau)}.
\end{align}
We seek to approximate $u(n,\beta,t)$ via an explicit time-stepping 
scheme.
Let $u_{k}(n,\beta)=u(n,\beta,t_{k})$, $N_{k}=N(t_{k})$,
$c_{k}=c(t_{k})$, $\xi_{k}=\xi(t_{k})$, and $t_{k+1}=t_{k}+\Delta t$.
Choosing a sufficiently small $\Delta t$ we compute: 

\vspace*{0.3cm}
{\bf Algorithm~1}
\begin{itemize}
\item{Given: $\xi_{0}=0$, $c_{0}>0$, $u_{0}(n,\beta)$, $\tilde{N_{0}}=0$,
and a sufficiently large $\Omega$.}
\item{For $k=0,1,2,\dots$, while $c_{k}\geq 0$, do:
\begin{align*}
&\xi_{k+1} = \xi_{k}+c_{k}\Delta t,
\\
&u_{k+1}(n,\beta) = e^{-\beta\xi_{k+1}}u_{0}(ne^{-\beta\xi_{k+1}},\beta),
\\
&N_{k+1} =  \iint\limits_{\Omega}nu_{k+1}(n,\beta)\,dn\,d\beta,
\\
&\tilde{N}_{k+1}=\tilde{N_{k}}+N_{k+1}\Delta t
\\
&c_{k+1}  = c_{0}\exp\left(- \gamma \tilde{N}_{k+1}\right)
\end{align*}
}
\end{itemize}
\begin{figure}[t]
\begin{center}
\hspace*{-0.5cm}
\includegraphics[width=10cm]{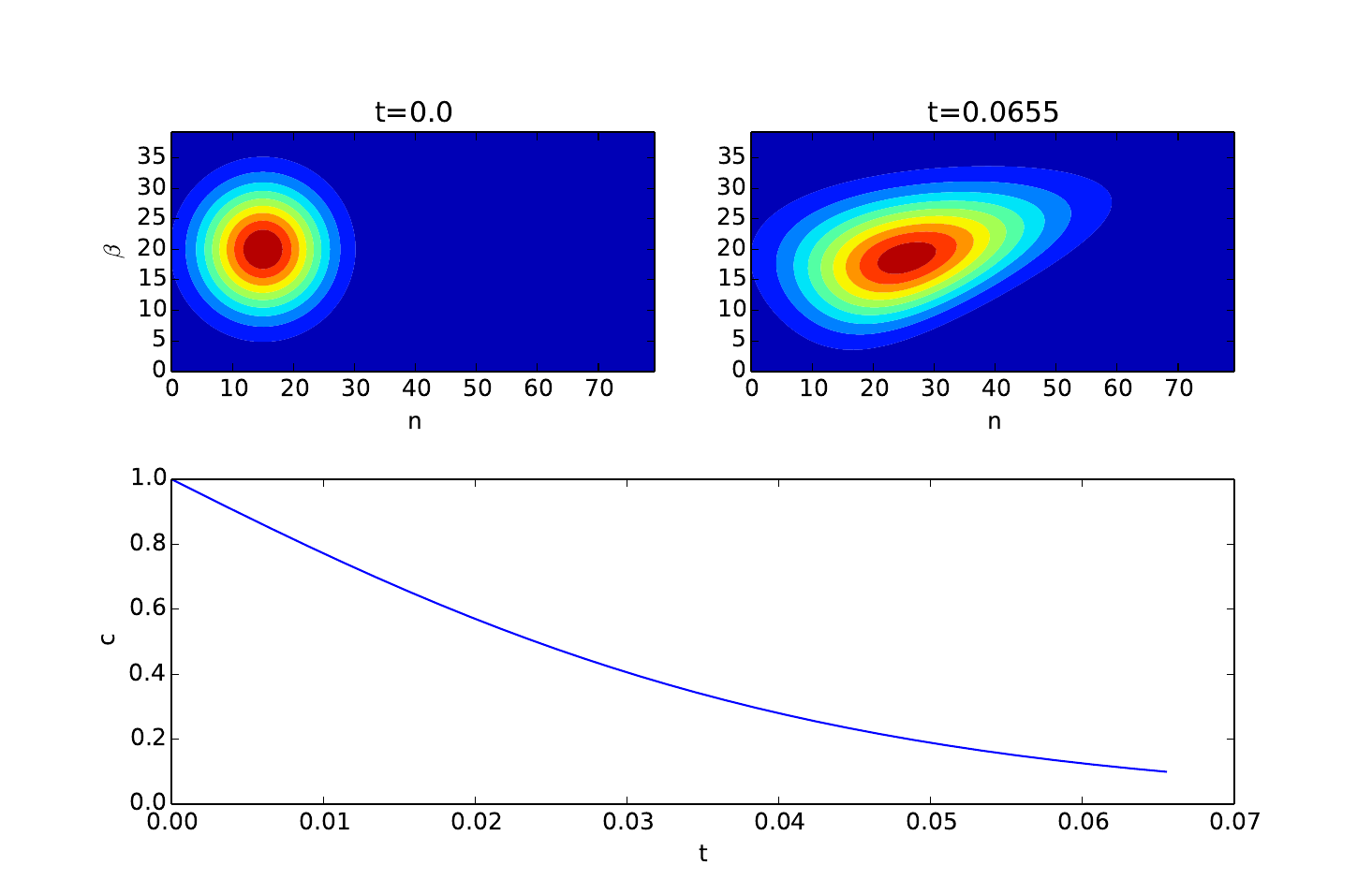}
\end{center}
\caption{
\label{fig:Fig3}
Growth of populations competing for a finite resource.
Top-left: initial distribution $u_{0}(n,\beta)$. Top-right:  
$u(n,\beta,t)$ for $t=0.0655$, at which point the amount of resource 
has dropped down to $c=0.1$ and the distribution 
practically does not change any more (arrested growth or suffocation due 
to the lack of resources).
Bottom: amount of the available resource $c(t)$ as a function of time.
}
\end{figure}

An example of the time evolution of $u(n,\beta,t)$ and the corresponding 
consumption of the resource $c(t)$ computed with the Algorithm~1 ($c_{0}=1$ and 
$\gamma=0.005$)
is shown in Fig.~\ref{fig:Fig3}. Due to the 
mathematical equivalence of the two solutions, 
the states (snapshots) of $u(n,\beta,t)$ given by (\ref{eq:uCompetingSolution}) are exactly 
equal to the states $u(n,\alpha,t)$ of the unmitigated exponential growth 
given by (\ref{eq:uSolution}), provided $\alpha=\beta$.
However, the evolution of $u(n,\beta,t)$ proceeds at a different rate 
and asymptotically slows down to a halt as $c(t)\rightarrow 0$. In fact,
it takes more than twice as much time to reach the same state 
as the last one of Fig.~\ref{fig:Fig1}.
The accuracy of the Algorithm~1
can be eventually improved using the predictor-corrector technique.

\subsection{Coexistence of competing species}
In the previous section all populations eventually stop growing (suffocate), 
since the resource $c$ is asymptotically exhausted for $t\rightarrow\infty$.
Here we assume that the resource is not only consumed, but is also generated at
a given fixed rate $\alpha$, i.e.,
\begin{align}
\label{eq:CoexN} 
\begin{split}
\frac{dc}{dt}&=\alpha-\gamma c \sum_{i}n_{i},\;\;\;c(0)=c_{0},\gamma>0,
\\
\frac{dn_{i}}{dt}&=\beta_{i}n_{i}(c-\kappa),\;\;\;n_{i}(0),\beta_{i},\kappa> 0,\;\;\;i=1,2,\dots
\end{split}
\end{align}
where $\kappa$ denotes the level of resource below which populations begin to decline 
(alternative interpretation is that for $c<\kappa$ the death rate becomes greater 
than the birth rate).
Different $\beta$'s -- reaction rates of species on the changes in the availability 
of resource (on the changes in the environment) -- could be viewed as 
different survival strategies. This problem is intrinsically interesting, because it
may result in the so-called equilibrium state, where species with different 
survival strategies coexist without exhausting the resource.

The analogue of the eq. (\ref{eq:nCompeting}) can be
obtained by integrating the first equation as
\begin{align}
\label{eq:cCoex}
c(t)=
c_{0}e^{-\gamma \tilde{N}(t)}+
\alpha e^{-\gamma \tilde{N}(t)}\int\limits_{0}^{t}e^{\gamma \tilde{N}(t')}\,dt',
\end{align}
and substituting this result into the second equation of (\ref{eq:CoexN}).
Note that $\tilde{N}$ is still given by (\ref{eq:tildeN}) and (\ref{eq:TotalPop}).
Keeping in mind that $c(t)$ is given by (\ref{eq:cCoex}), we write the phase-space 
problem as
\begin{align}
\label{eq:uCoexisting}
\frac{\partial u}{\partial t}+
(c-\kappa)\beta n\frac{\partial u}{\partial n} + (c-\kappa)\beta u = 0, \;\;\;u(n,\beta,0)=u_{0}(n,\beta),
\end{align}
Mathematically this problem is more challenging than the one of the previous section,
since the phase-space velocity $(c-\kappa)\beta n$ may change its sign, albeit for 
all $n$ and $\beta$ at the same time. Moreover, one cannot immediately 
conclude that this velocity will tend to some well-defined limit as $t\rightarrow\infty$.

The following Lemma~\ref{lem1} establishes the necessary condition for an instantaneous 
{\it steady-state} at time $t$,
which we define as $u(n,\beta,t)$, such that $\partial u/\partial t = 0$ at $t$.
This is the moment when the phase-space velocity becomes zero everywhere in $\Omega$.

\begin{lemma}
\label{lem1}
A function $\tilde{u}(n,\beta,t)\neq 0$, integrable on $\Omega=\{(n,\beta)\vert\,
0\leq n<\infty,\,0\leq\beta<\infty\}$ for any $t\geq 0$, 
is a steady-state solution of the 
problem (\ref{eq:uCoexisting}) at time $t$ if and only if this function 
yields $c(t)=\kappa$ with $c(t)$ given by (\ref{eq:cCoex}), (\ref{eq:TotalPop}).
\end{lemma}
\begin{proof} 
The if part is obvious. To prove the rest, assume that there exists a steady-state 
solution $\tilde{u}(n,\beta,t)\neq 0$ integrable on $[0,\infty)\times[0,\infty)$ for
$c(t)\neq\kappa$. Hence, it must satisfy the following linear equation:
\begin{align}
\label{eq:Lem1Eq1}
n\frac{\partial\tilde{u}}{\partial n}+\tilde{u} = 0,
\end{align}
i.e., it should have the form $\tilde{u}=C/n$, where $C\neq 0$ is a constant. 
Such functions, however, are not integrable on 
$[0,\infty)\times[0,\infty)$, and we arrive at a contradiction.
\end{proof}

A steady state defined as above is, in general, not stable. 
Indeed, although $\partial u/\partial t=0$ whenever $c(t)=\kappa$,
the resource function $c(t)$ may continue to change in time, 
thus deviating away from its steady-state level $\kappa$, which,
according to Lemma~1, will necessarily `restart' the evolution of $u$. 
The time derivative of $c(t)$ is found from (\ref{eq:cCoex}) to be
$c'=\alpha-\gamma c N$ and, understandably, coincides with the 
resource equation in (\ref{eq:CoexN}).
In a steady state we have $c'(t)=\alpha-\gamma \kappa N(t)$, which
may be non-zero. To prevent $u$ from leaving the steady state, one needs $c'(t)=0$
as well, i.e., $N(t)=\alpha/(\gamma\kappa)$. 

\begin{definition}
\label{def1}
A function $u(n,\beta,t)$ represents an {\it equilibrium} state 
if it simultaneously satisfies the following conditions:
\begin{align}
\label{eq:Def1}
c(t)&=\kappa,
\\
\label{eq:Def2}
N(t)&=\frac{\alpha}{\gamma\kappa},
\end{align}
where $c(t)$ and $N(t)$ depend on $u$ as in (\ref{eq:cCoex}) and (\ref{eq:N}).
\end{definition}

We shall assume the existence and uniqueness
of the solution $u(n,\beta,t)$ to the problem (\ref{eq:cCoex})--(\ref{eq:uCoexisting})
for any $t\geq 0$. It should be possible to arrive at 
the corresponding proof along the lines of Collet and Goudon (2000) who considered
very similar problems. 
Here we focus on the stability of the equilibrium state.

\begin{theorem}
\label{Thm1}
Let $\alpha,\gamma,\kappa>0$, and let $u(n,\beta,t)$ be a solution 
of (\ref{eq:uCoexisting}) on $\Omega=\{(n,\beta)\vert\,
0\leq n<\infty,\,0\leq\beta<\infty\}$ for any $t\geq 0$, 
corresponding to the initial state $u_{0}(n,\beta)$,
such that $u_{0}(n,\beta)=0$ for $\beta > B$. Further, let
\begin{align}
\label{eq:Thm1Ass1}
\lim_{n\rightarrow\infty}n^{2} u(n,\beta,t) = 0,
\;\;\;\beta\in[0,B],\;\;\;t\geq 0.
\end{align} 
Then, the equilibrium state is asymptotically stable.
\end{theorem}
\begin{proof}
To analyze the stability of equilibrium we need 
another direct relation between $c(t)$ and $N(t)$ in addition to (\ref{eq:cCoex}),
preferably, less complicated than the PDE (\ref{eq:uCoexisting}). 
For this purpose we rewrite the time derivative of $N(t)$ as:
\begin{align}
\label{eq:dNdtDerivation}
\begin{split}
&N'=\iint_{\Omega}n\frac{\partial u}{\partial t}\,d\Omega=
-(c-\kappa)\iint_{\Omega}\beta n\frac{\partial}{\partial n}
\left(n u\right)\,d\Omega 
\\
&= -(c-\kappa)\iint_{\Omega}
\beta\frac{\partial}{\partial n}\left(n^{2} u\right)\,d\Omega +
(c-\kappa)\iint_{\Omega}\beta n u \,d\Omega
\\
&= -(c-\kappa)\overline{\beta}B\lim_{\nu\rightarrow\infty}\int_{0}^{\nu}
\frac{\partial}{\partial n}\left(n^{2} u(n,\overline{\beta},t)\right)\,dn +
(c-\kappa)\beta_{\rm c}\iint_{\Omega}n u \,d\Omega
\\
& = -(c-\kappa)\overline{\beta}B\lim_{\nu\rightarrow\infty}
\left[\nu^{2} u(\nu,\overline{\beta},t)\right] + (c-\kappa)\beta_{\rm c}N,
\end{split}
\end{align}
where we have employed the Mean Value Theorem with 
$\overline{\beta},\beta_{\rm c}\in [0,B]$.
Hence, using (\ref{eq:Thm1Ass1}) we arrive at the coupled quasi-linear system:
\begin{align}
\label{eq:dNdt}
\begin{split}
c'& = \alpha - \gamma c N,
\\
N'& =(c-\kappa)\beta_{\rm c}N,
\end{split}
\end{align}
with the equilibrium (\ref{eq:Def1})--(\ref{eq:Def2}) as its critical point. 
The Jacobian at equilibrium is given by:
\begin{align}
\label{eq:Jacobian}
\begin{bmatrix}
-\alpha/\kappa & -\gamma\kappa\\
\alpha\beta_{\rm c}/(\gamma\kappa) & 0
\end{bmatrix},
\end{align} 
and its eigenvalues are:
\begin{align}
\label{eq:EigsCN}
\lambda_{1,2}=-\frac{\alpha}{2\kappa}
\pm \frac{\sqrt{\alpha(\alpha-4\beta_{\rm c}\kappa^{2})}}{2\kappa}.
\end{align}
The real part of the eigenvalues is negative for all $\alpha,\kappa,\beta_{\rm c}>0$. 
Since, in equilibrium both $\partial u/\partial t =0$ and $c'=0$, 
stability of equilibrium for $c(t)$ and $N(t)$ implies
stability of equilibrium for $u(n,\beta,t)$. 
\end{proof}

Thus, we have established that, as soon as the 
distribution $u(n,\beta,t)$ evolves into a state, for which
the scalars $c(t)$ and $N(t)$ are close enough to their equilibrium values
(\ref{eq:Def1})--(\ref{eq:Def2}), the system (\ref{eq:cCoex})--(\ref{eq:uCoexisting})
will converge to the equilibrium.
The possibility of complex eigenvalues for $\alpha < 4B\kappa^{2}$ 
shows that $c$ and $N$ may spiral towards the equilibrium point, i.e.,
oscillate in time. The condition (\ref{eq:Thm1Ass1}) of the Theorem~\ref{Thm1},
albeit natural (it means the absence of infinitely large populations at any given $t\geq 0$), 
is not really necessary. Everything depends on how close the system is
to its equilibrium state at $t=0$. It is sufficient that the convergence to equilibrium
happens faster in time than the blow up of potentially unbounded populations.

To investigate the behavior of the distribution function for coexisting species 
we employ the formal implicit solution of (\ref{eq:uCoexisting}), which 
has the same form as (\ref{eq:uCompetingSolution}) with
\begin{align}
\label{eq:xiCoexisting}
\xi(t)=\int\limits_{0}^{t}(c(t')-\kappa)\,dt',
\end{align}
where $c(t)$ is given by (\ref{eq:cCoex}). The analogue of the Algorithm~1
in the present case may be formulated as follows:

\vspace*{0.3cm}
{\bf Algorithm~2}
\begin{itemize}
\item{Given: $\xi_{0}=0$, $c_{0}>0$, $u_{0}(n,\beta)$, $\tilde{N_{0}}=0$, $R_{0}=0$,
and a sufficiently large $\Omega$.}
\item{For $k=0,1,2,\dots$, while $\vert c_{k}-\kappa\vert>0$ and 
$\vert N_{k}-\alpha/(\gamma\kappa)\vert>0$, do:
\begin{align*}
&\xi_{k+1} = \xi_{k}+\Delta t (c_{k}-\kappa),
\\
&u_{k+1}(n,\beta) = e^{-\beta\xi_{k+1}}u_{0}(ne^{-\beta\xi_{k+1}},\beta),
\\
&N_{k+1} =  \iint\limits_{\Omega}nu_{k+1}(n,\beta)\,dn\,d\beta,
\\
&\tilde{N}_{k+1}=\tilde{N_{k}}+\Delta t N_{k+1}
\\
&R_{k+1}=R_{k}+\Delta t\exp\left(\gamma \tilde{N}_{k+1}\right)
\\
&c_{k+1}  = \left(c_{0}+\alpha R_{k+1}\right)\exp\left(- \gamma \tilde{N}_{k+1}\right)
\end{align*}
}
\end{itemize}
As a measure of approach to equilibrium we use the differences 
$\vert c_{k}-\kappa\vert$ and $\vert N_{k}-\alpha/(\gamma\kappa)\vert$.
\begin{figure}[t]
\begin{center}
\hspace*{-0.5cm}
\includegraphics[width=12cm]{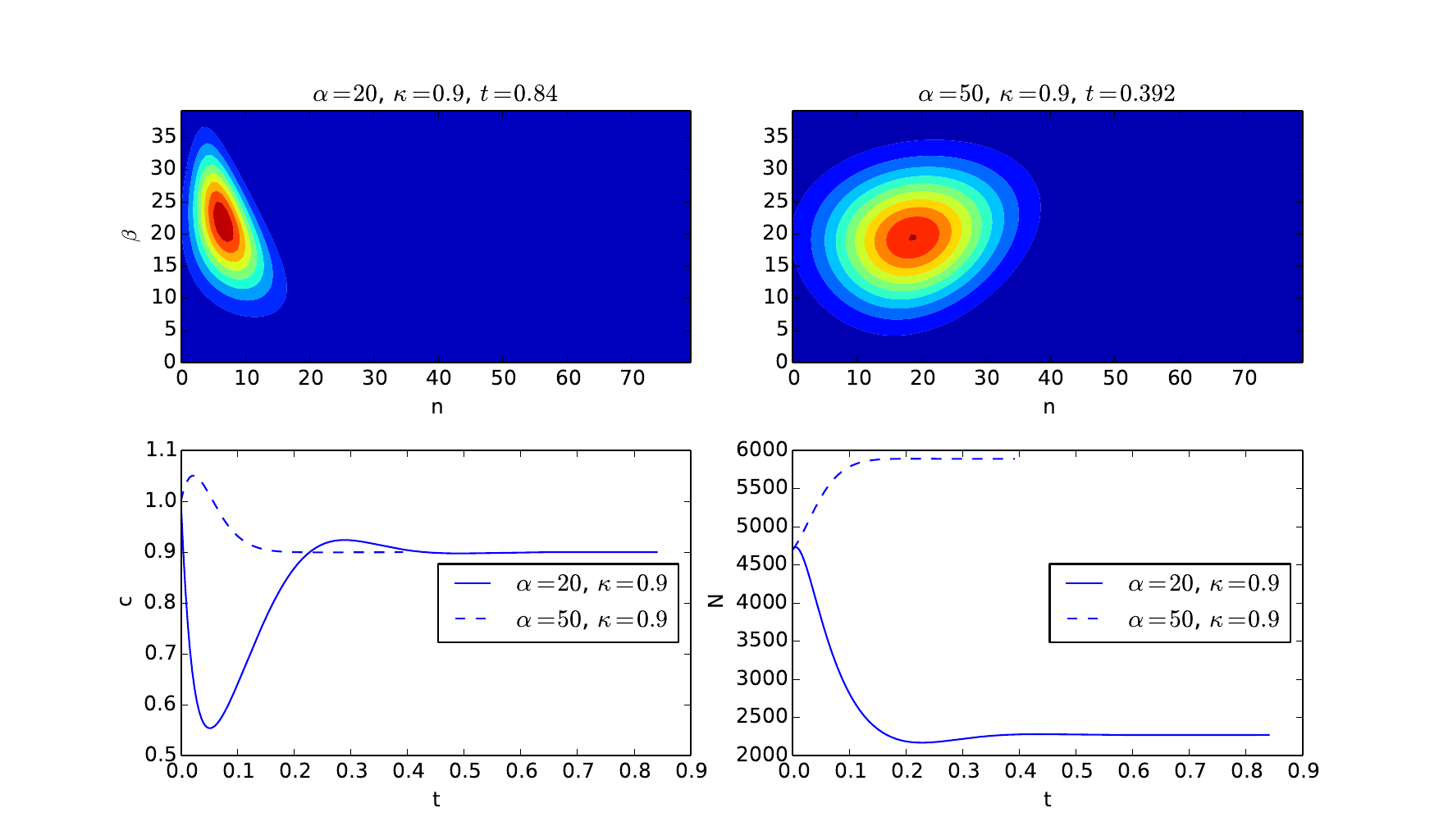}
\end{center}
\caption{
\label{fig:Fig4a}
Coexistence of competing species.
Top: quasi-equilibrium distributions $u(n,\beta,t)$
achieved starting with the same initial distribution as in Fig.~\ref{fig:Fig3}
for environments with different resource generation rates 
(left: $\alpha = 20$, right: $\alpha = 50$).
The distributions are shown at $t=0.84$ and $t=0.392$, respectively.
Bottom-left: amount of the available resource $c(t)$ as a function of time. Bottom-right:
the time evolution of the total size of all populations $N(t)$.
}
\end{figure}

Figure~\ref{fig:Fig4a} illustrates two cases of the evolution of the distribution 
function $u(n,\beta,t)$ for competing species with the parameters of 
the problem set as $c_{0}=1.0$, $\gamma=0.01$, $\kappa = 0.9$, 
$\Delta t = 0.001$, and the same initial distribution as in 
Figure.\ref{fig:Fig3}. 
 
As expected, in all our numerical experiments 
the distribution functions (virtually) stop changing as soon as 
the resource and the total size approach their equilibrium values. 
For example, the two images in the top 
of Figure~\ref{fig:Fig4a} and the solid and dashed curves 
in the bottom plots show the stabilization of the distribution in environments
with different resource generation rates: $\alpha = 20$ corresponds to the 
top-left image and solid curves, $\alpha = 50$ corresponds to the top-right 
image and dashed curves. 
As predicted, after a few oscillations the resource and the total size 
tend to stabilize. Although, the total size converges 
to $N = \alpha/(\gamma\kappa)$ at a much slower rate than the resource
converges to $c=\kappa$.

\begin{figure}[t]
\begin{center}
\hspace*{-0.5cm}
\includegraphics[width=12cm]{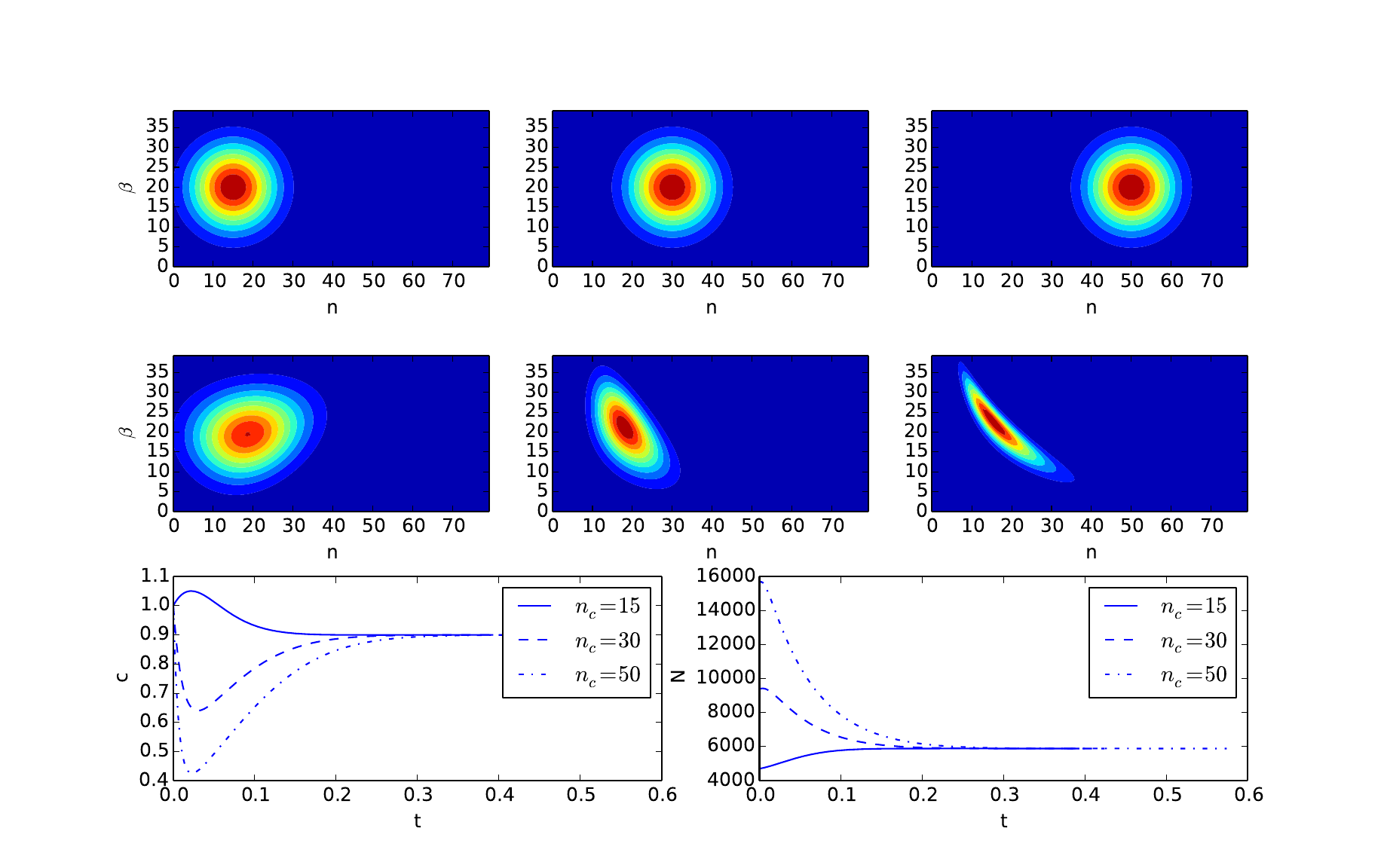}
\end{center}
\caption{
\label{fig:Fig4b}
Dependence of the quasi-equilibrium for coexisting species on the initial state.
Top row: initial distributions $u_{0}(n,\beta)$ with centers at $n_{c}=15$,
$n_{c}=30$, and $n_{c}=50$.
Middle row: quasi-equilibrium distributions obtained from the above initial distributions.
Bottom-row: amount of the available resource $c(t)$ and the total size of all 
populations $N(t)$ as functions of time.
}
\end{figure}

Although the evolution towards equilibrium appears to be 
a stable process, the equilibrium state itself is not unique (even with 
fixed parameters $c_{0}$, $\gamma$, $\alpha$,
and $\kappa$) and strongly depends on the initial state $u_{0}$. 
Figure~\ref{fig:Fig4b} illustrates this phenomenon for $c_{0}=1$,
$\gamma=0.01$, $\alpha=50$, and $\kappa=0.9$, and three different initial states
(Gaussian distributions with their centers at $\beta_{\rm c}=20$,
and $n_{\rm c}=15$, $30$, and $50$, respectively).
The resource $c$ and the total size of populations $N$ approach the same 
values in all three cases. However, the quasi-equilibrium distributions are
very different. 

In general, however, a stable coexistence of species whose survival strategies 
span a continuum of possibilities is not possible (Gyllenberg and Mesz{\'e}na ~2005). This
result can be confirmed in the present formulation as well if one allows for 
different carrying capacities. Namely, consider a generalized version of the
problem (\ref{eq:CoexN}):
\begin{align}
\label{eq:GeneralizedCoexN} 
\begin{split}
\frac{dc}{dt}&=\alpha- c \sum_{i}\gamma_{i}n_{i},\;\;\;c(0)=c_{0},\;\;\gamma_{i}>0,\;\;\;i=1,2,\dots,
\\
\frac{dn_{i}}{dt}&=\beta_{i}n_{i}(c-\kappa_{i}),\;\;\;n_{i}(0),\beta_{i},\kappa_{i}> 0,\;\;\;i=1,2,\dots,
\end{split}
\end{align}
where the rate of resource consumption $\gamma_{i}$ as well as the lowest growth-sustaining 
level of resource $\kappa_{i}$ depend on the species. Multiplying the second equation in 
(\ref{eq:GeneralizedCoexN}) with $\gamma_{i}$ and introducing the new size variable
$m_{i}=\gamma_{i}n_{i}$ we arrive at the following problem:
\begin{align}
\label{eq:GeneralizedCoexM} 
\begin{split}
\frac{dc}{dt}&=\alpha- c \sum_{i}m_{i},\;\;\;c(0)=c_{0},
\\
\frac{dm_{i}}{dt}&=\beta_{i}m_{i}(c-\kappa_{i}),\;\;\;m_{i}(0),\beta_{i},\kappa_{i}> 0,\;\;\;i=1,2,\dots,
\end{split}
\end{align}
which has almost the same form as (\ref{eq:CoexN}), except for the species-dependent 
$\kappa_{i}$. The corresponding distribution function $u(m,\beta,\kappa,t)$ 
satisfies the same phase-space equations (\ref{eq:cCoex})--(\ref{eq:uCoexisting}) 
with $n$ replaced by $m$, $\gamma=1$, and with $N(t)$, $\tilde{N}(t)$
replaced by $M(t)$, $\tilde{M}(t)$, where the latter are defined as:
\begin{align}
\label{eq:M}
M(t)=\iiint\limits_{\Omega} m u(m,\beta,\kappa,t)\,d\Omega,
\;\;\;
\tilde{M}(t)=\int\limits_{0}^{t}M(t')\,dt'.
\end{align}
Obviously, Lemma~\ref{lem1} prohibits any equilibrium for this 
problem, since there does not exist a single value of the resource $c(t)$ 
that would satisfy the necessary condition $c(t)=\kappa$ over a whole
range of $\kappa$'s. Thus, the distribution function will keep changing 
in time.

\subsection{Randomized migration}
Let a single metapopulation be nonuniformly distributed over some spatial domain at time $t=0$.
The whole habitat may be divided either naturally or virtually into a large number of 
elementary cells (sub-habitats) and a question can be posed about the distribution of cell populations 
over their size at time $t$.
The growth of these elementary cell populations could be described by one of the models from the
preceding sections with growth parameters varying from cell to cell. If all individual members
stay within their original cells or the whole metapopulation uniformly translates in space
the phase-space methods developed earlier 
and the resulting distribution functions, obviously, apply without change. 

Allowing for the migration of individuals between the cells
usually requires a PDE-based space-time dynamic law (e.g. Fischer equation). 
Such a law cannot be directly
incorporated into the phase-space current, since a low-dimensional phase-space approach 
completely neglects the spatial ordering 
of the cells. Nevertheless, certain types of migration mechanisms can
be described by ODE systems with low-dimensional phase spaces.

Consider a stationary metapopulation, where the total number of individuals 
and the total number of sub-habitats (that play the role of populations in this case) 
do not change. Suppose that all individuals from time to time (with frequency $\beta$) 
decide to emigrate from 
its native population so that each population is subject to the emigration rate $\beta n$,
proportional to its size. Depending on the choice strategy of migrants as far as their
new population is concerned one arrives at different ODE systems.

For example, let each emigrant choose a new population (habitat) completely at random.
If viewed spatially, this is, obviously, a variant of super-diffusion with a random step 
size, which is a non-trivial problem if one considers it directly in space and time.
On a sufficiently large time scale, there will be the steady immigration 
rate $\gamma$ in all populations (since the total emigration rate is constant). 
The emigration and immigration coefficients $\beta$ and $\gamma$ must be balanced to
guarantee the steady state of the metapopulation:
\begin{align}
\label{eq:BalanceRandom}
\frac{d}{dt}\sum\limits_{i=1}^{P}n_{i}
=\sum\limits_{i=1}^{P}(\gamma-\beta n_{i}) = 0.
\end{align}
Let $\sum_{i}n_{i}=N={\rm const}$, then $\gamma = \beta \overline{n}$, where $\overline{n}=N/P$, 
and $P$ is the number of habitats. The dynamic law becomes
\begin{align}
\label{eq:EmigationRandom}
\frac{dn_{i}}{dt}=\beta\left(\overline{n} - n_{i}\right),\;\;\;i=1,2,\dots P.
\end{align}
Notice that $\overline{n}=N/P$ is the total number of individuals divided by the number of
habitats, which would be exactly the size of populations in all habitats if all individuals
were evenly spread throughout. From (\ref{eq:EmigationRandom}) it follows that a population
stops changing once it reaches this equilibrium size $\overline{n}$.

The distribution function $u(n,t)$ is then found by solving the following
Cauchy problem:
\begin{align}
\label{eq:ContinuityEmigrationRandom}
\frac{\partial u}{\partial t} + 
\beta\left(\overline{n} - n\right)\frac{\partial u}{\partial n} 
-\beta u = 0,
\;\;\; u(n,0)=u_{0}(n),
\end{align}
and is given by
\begin{align}
\label{eq:uEmigrationRandom}
u(n,t)=
\begin{cases}
e^{\beta t}u_{0}\left(\overline{n}-(\overline{n}-n)e^{\beta t}\right), & n<\overline{n},
\\
e^{\beta t}u_{0}\left(\overline{n}\right), & n=\overline{n},
\\
e^{\beta t}u_{0}\left(\overline{n}+(n-\overline{n})e^{\beta t}\right), & n>\overline{n}.
\end{cases}
\end{align}
\begin{figure}[t]
\begin{center}
\hspace*{-0.5cm}
\includegraphics[width=13cm]{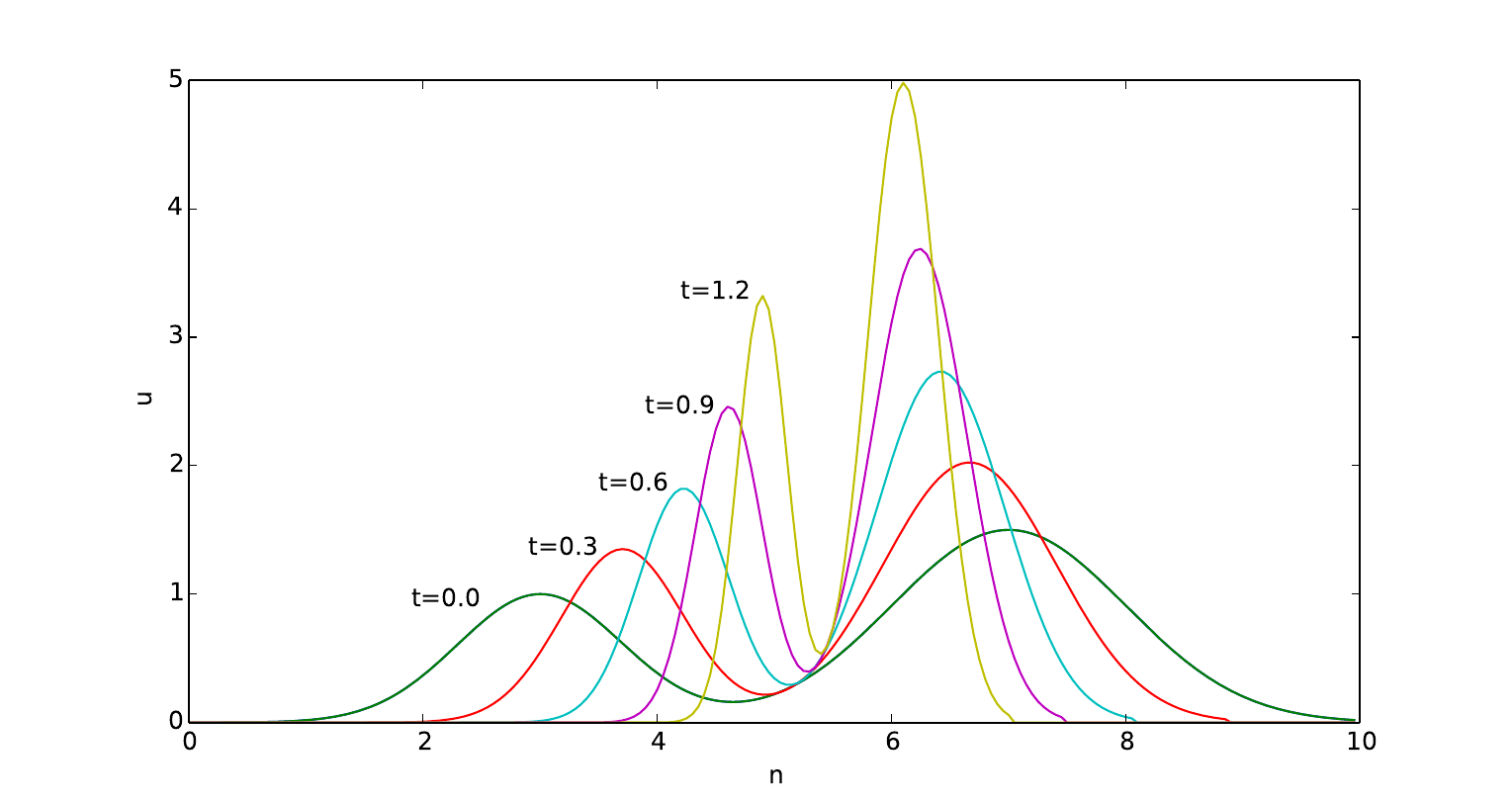}
\end{center}
\caption{
\label{fig:Fig5}
Time evolution of the distribution function $u(n,t)$ for populations with (uniform) migration
given by the equation (\ref{eq:uEmigrationRandom}) with $\beta = 1.0$ and $\overline{n}=5.7$.
}
\end{figure}
The middle branch of this solution shows that the number (density) of populations with $n=\overline{n}$
growth exponentially with time. As can be seen in Fig.~\ref{fig:Fig5} the upper and lower branches
of the solution (\ref{eq:uEmigrationRandom}) 
represent the initial distribution that shrinks, respectively, from the left and
from the right, towards $\overline{n}$. This behaviour is indicative of a mollifier of the Dirac 
delta function centred at $\overline{n}$, which, indeed, is the limit of 
$u(n,t)$ as $t\rightarrow\infty$, since all populations will eventually have the same size 
$\overline{n}$. An interesting feature of (\ref{eq:uEmigrationRandom})
is the preservation of the main features of the initial $u_{0}(n)$ over time, albeit in a scaled (shrank) form,
e.g., the double peaked shape shown in Fig.~\ref{fig:Fig5}. 


\subsection{Biased migration}
Consider populations with potentially unlimited growth and constant emigration frequency.
Let the choice of the new population by the migrants be biased towards
populations with higher growth rates. The ODE describing this situation is:
\begin{align}
\label{eq:EmigrationBiased}
\frac{d n_{i}}{dt} = \alpha_{i}n_{i}-\beta n_{i} + \gamma \alpha_{i},
\;\;\;i=1,\dots P,
\end{align}
where $\alpha_{i}$ is the growth coefficient, $\beta$ is the emigration 
frequency, and $\gamma$ is the immigration coefficient. Since migration does not
change the total number of individuals (only the growth does), the following balance equation 
should be satisfied at all times:
\begin{align}
\label{eq:BalanceBiased}
\beta \sum\limits_{i=1}^{P} n_{i} = \gamma \sum\limits_{i=1}^{P} \alpha_{i},
\end{align}
reducing equation (\ref{eq:EmigrationBiased}) to
\begin{align}
\label{eq:EmigrationBiased1}
\frac{d n_{i}}{dt} = (\alpha_{i}-\beta) n_{i} + \alpha_{i}\beta
\frac{\sum n_{j}}{\sum \alpha_{j}},
\;\;\;i,j=1,\dots P.
\end{align}
The corresponding nonzero component of the phase-space current is
\begin{align}
\label{eq:JBiased}
J_{n}=(\alpha-\beta)nu + \alpha\beta Ru,
\end{align}
where
\begin{align}
\label{eq:NBiased}
\begin{split}
R&=\frac{N}{A},
\\
N=\iint n u(n,\alpha,t)\,d\alpha\,dn,
&\;\;\;
A=\iint \alpha u(n,\alpha,t)\,d\alpha\,dn,
\end{split}
\end{align}
and, since the $\alpha$-component of the phase-space current is zero, $A$ will stay constant 
in time.
Thus, we arrive at the following weakly nonlinear continuity equation for the distribution 
$u(n,\alpha,t)$:
\begin{align}
\label{eq:ContinuityBiased}
\frac{\partial u}{\partial t} + (\alpha\beta R+(\alpha-\beta)n)\frac{\partial u}{\partial n}
+(\alpha-\beta)u = 0,\;\;\;u(n,\alpha,0)=u_{0}(n,\alpha).
\end{align}
Assuming that $R(t)$ is a given function of time we use the method of characteristics
to obtain the implicit solution of this problem as
\begin{align}
\label{eq:uBiased}
\begin{split}
u(n,\alpha,t)&=e^{-(\alpha-\beta)t}
u_{0}\left(ne^{-(\alpha-\beta)t}-\alpha\beta\xi(\alpha,t),\alpha\right),
\\
\xi(\alpha,t)&=\int\limits_{0}^{t}R(t')e^{-(\alpha-\beta)t'}\,dt'.
\end{split}
\end{align}
However, $R(t)$ and, hence, $\xi(\alpha,t)$ are functions of $u$ as well, see eq.~(\ref{eq:NBiased}). 
Thus, we employ an iterative algorithm similar to the Algorithms~1 and 2.

\vspace*{0.3cm}
{\bf Algorithm~3}
\begin{itemize}
\item{Given: $\xi_{0}(\alpha)=0$, $u_{0}(n,\beta)$, $N_{0}$, $R_{0}$,
and a sufficiently large $\Omega$.}
\item{For $k=0,1,2,\dots$, do:
\begin{align*}
&\xi_{k+1}(\alpha) = \xi_{k}(\alpha)+R_{k}e^{-(\alpha-\beta)t_{k}}\Delta t,
\\
&u_{k+1}(n,\alpha)=e^{-(\alpha-\beta)t_{k+1}}
u_{0}\left(ne^{-(\alpha-\beta)t_{k+1}}-\alpha\beta\xi_{k+1}(\alpha),\alpha\right),
\\
&R_{k+1} =  \frac{1}{A}\iint\limits_{\Omega}nu_{k+1}(n,\alpha)\,dn\,d\alpha,
\end{align*}
}
\end{itemize}
\begin{figure}[t]
\begin{center}
\hspace*{-0.5cm}
\includegraphics[width=13cm]{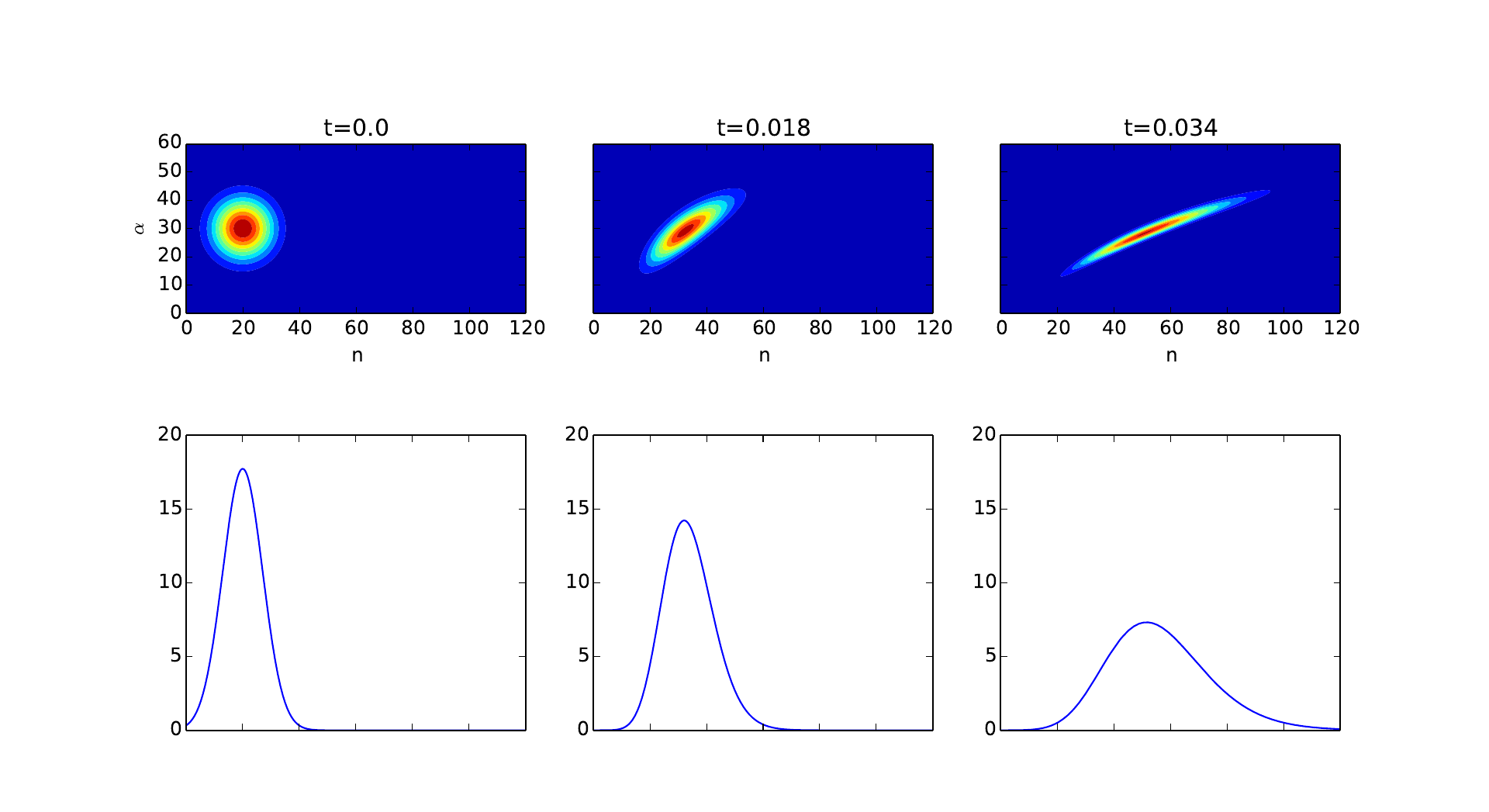}
\end{center}
\caption{
\label{fig:Fig6}
Time evolution of the distribution function $u(n,\alpha,t)$ for populations 
with biased migration, see eq.~(\ref{eq:uBiased}). 
Top row -- snapshots at $t=0.0, 0.018$, and $0.034$; Bottom 
row -- corresponding numerically integrated $\rho(n,t)$ showing the distribution of
populations over their sizes.
}
\end{figure}

A quick look at the implicit solution (\ref{eq:uBiased}) reveals similarity with
the case of exponential growth (\ref{eq:uUnlimited}). Indeed, as expected, the two solutions are 
identical for $\beta=0$, remain similar at low values 
of the emigration frequency, and become very different for larger values.
One of the predictions of the metapopulation theory (Levin~1969; Eriksson et~al.~2014)
is that migration may prevent extinction. In the present case, as follows from
(\ref{eq:EmigrationBiased1}), none of the populations can disappear as long as there 
is at least one nonzero population somewhere. Yet, it does not mean that all
populations are always growing in size. In the early stages some populations may
decrease. From (\ref{eq:EmigrationBiased1}) we deduce that the condition on the 
growth of a population is
\begin{align}
\label{eq:GrowthCondition}
\frac{\alpha_{i}}{\sum_{j}\alpha_{j}}
>\left(1-\frac{\alpha_{i}}{\beta} \right)\frac{n_{i}}{\sum_{j}n_{j}},\;\;\;\beta>0.
\end{align}
Obviously, this condition is satisfied for all populations with $\alpha_{i}>\beta$. 
However, for a population with $\alpha_{i}<\beta$ it is easy to imagine a
sufficiently large initial state $n_{i}(0)$ such that this condition is violated 
and the population size will decrease. Nevertheless, when such a decreasing $n_{i}$ 
becomes sufficiently small, the condition will become satisfied again and the population 
will resume its growth.

An example of the evolution of distribution function for populations with biased 
migration computed with the Algorithm~3 is shown in Fig.~\ref{fig:Fig6}.
The initial distribution (the leftmost image and plot) is a Gaussian centred at
$n_{\rm c}=20$ and $\alpha_{\rm c}=30$. The emigration frequency is set at 
$\beta=50$ in this simulation, which is higher than most of the growth 
coefficients $\alpha$ in the initial distribution. The first striking 
conclusion is that despite the immigration bias we do not observe an explosive growth 
of some particular (chosen) populations at the expense of others. Apparently, any such 
growth is efficiently mitigated by a higher emigration rate. 

The early-time transient decrease of some populations mentioned above 
(those with large initial $n$'s) 
causes the distribution function to contract slightly in the horizontal $n$-dimension,
since the right side of the distribution moves to the left, while the left side moves 
to the right. This happens in the very early stages of the evolution (not shown). 
The same effect prevents the distribution from widening in the horizontal direction 
in the course of the evolution. In fact, as the Figure~\ref{fig:Fig6} shows, 
the distribution only keeps contracting.

Thus, the main consequence of the large emigration frequency $\beta$ for the 
late-time evolution appears to be the uniformity of size for populations with 
equal $\alpha$'s and the uniformity of $\alpha$'s for populations of the same size.
Also, for large $t$'s the populations show a very strong
correlation (almost direct proportionality) between the population size 
and its growth coefficient $\alpha$ (see the top-right image of Fig.~\ref{fig:Fig6}).

\subsection{Growth rates depending on the distribution function}
The present approach is indispensable for 
dynamic models where the rate of growth of populations is a 
functional $R(u,n)$ of both the size $n$ of the local population and the distribution $u$ 
of populations over the phase space. 
The phase-space current vector will have as many nonzero 
components as there are time varying parameters in $R$. 
If all parameters,
except $n_{i}$, are time invariant, then the problem is described by the
following equations:
\begin{align}
\label{eq:GeneralDynamicDiscrete}
\frac{d n_{i}}{d t} & = R(u,n_{i}),\;\;\;i=1,2,\dots P,
\\
\label{eq:GeneralContinuity}
\frac{\partial u}{\partial t} & = 
- \frac{\partial}{\partial n}\left(uR(u,n)\right).
\end{align}
Obviously, it is now impossible to solve the dynamic equations 
(\ref{eq:GeneralDynamicDiscrete}) without first solving the phase-space 
problem (\ref{eq:GeneralContinuity}). 

Consider populations with unlimited growth discussed earlier. 
Suppose that one is able to influence 
the rate of growth of each population (e.g., by judicially watering and fertilizing each plant) 
proportionally to its current `weight' in the phase space. 
Such a problem is described by the following equations:
\begin{align}
\label{eq:DynamicNU}
\frac{d n_{i}}{d t} & = 
\alpha_{i}n_{i}
\int\limits_{\alpha_{i}-\Delta\alpha}^{\alpha_{i}+\Delta\alpha}
\int\limits_{n_{i}-\Delta n}^{n_{i}+\Delta n}u(n,\alpha,t)\,dn\,d\alpha
\approx \alpha_{i}n_{i}u(n_{i},\alpha_{i},t)\delta,\;\;\;i=1,2,\dots P,
\\
\label{eq:ContinuityNU}
\frac{\partial u}{\partial t} & +  
2 \alpha \delta\, n u \frac{\partial u}{\partial n} 
+ \alpha \delta\, u^{2} = 0, \;\;\;u(n,\alpha,0)=u_{0}(n,\alpha),
\end{align}
where $\delta=4\Delta\alpha\, \Delta n$, i.e.,
we use the mid-point approximation for the integral in (\ref{eq:DynamicNU}).
The resulting approximate phase-space problem (\ref{eq:ContinuityNU}) is a 
variant of the Burgers equation (see e.g. Smoller~1994) whose characteristic equations are:
\begin{align}
\label{eq:CharacteristicNU}
\frac{dt}{ds}=1,\;\;\;
\frac{dn}{ds}=2\alpha\delta\, n u,\;\;\;
\frac{du}{ds}=-\alpha\delta\, u^{2},
\end{align}
or
\begin{align}
\label{eq:CharacteristicNU2}
dn = 2\alpha\delta\, n u\, dt,\;\;\;
-\alpha\delta\, dt=\frac{du}{u^{2}},\;\;\;
\frac{dn}{2n}=-\frac{du}{u}.
\end{align}
The solutions of the last two are:
\begin{align}
\label{eq:CharacteristicNU3}
\alpha\delta\, t -\frac{1}{u} = C_{2},\;\;\;
u\sqrt{n}=C_{3},
\end{align}
which allows rewriting and solving the first equation in (\ref{eq:CharacteristicNU2}) as
\begin{align}
\label{eq:CharacteristicNU4}
\frac{dn}{\sqrt{n}} = 2\alpha\delta C_{3}\, dt,\;\;\;
2\sqrt{n}(1-\alpha\delta\, t u)=C_{1}.
\end{align}
Combining this result with the first equation in (\ref{eq:CharacteristicNU3}) we get
\begin{align}
\label{eq:CharacteristicNU5}
u=\frac{1}{\alpha\delta\, t - f(p)},\;\;\;p=2\sqrt{n}(1-\alpha\delta\, t u).
\end{align}
The initial condition requires that
\begin{align}
\label{eq:CharacteristicNU6}
\left.f(p)\right\vert_{t=0}=-\frac{1}{u_{0}(n,\alpha)}.
\end{align}
Hence, we choose
\begin{align}
\label{eq:CharacteristicNU7}
f(p)=-\frac{1}{u_{0}(p^{2}/4,\alpha)}=-\frac{1}{u_{0}(n(\alpha\delta\, t u-1)^{2},\alpha)},
\end{align}
and arrive at the following equation that implicitly defines the distribution $u(n,\alpha,t)$:
\begin{align}
\label{eq:uNU1}
u=\frac{u_{0}(n(\alpha\delta\,t\,u-1)^{2},\alpha)}
{\alpha\delta\,t\,u_{0}(n(\alpha\delta\,t\,u-1)^{2},\alpha) + 1},
\end{align}
or 
\begin{align}
\label{eq:uNU2}
u=(1-\alpha\delta\,t\,u)\,u_{0}(n(1-\alpha\delta\,t\,u)^{2},\alpha).
\end{align}
These representations are only valid in the smooth regime, i.e., prior to the 
development of shock.

A fixed-point iterative algorithm that computes the distribution at a given time $t$ 
may be formulated as follows (one can use any suitable norm for the computation 
of the current mismatch $\epsilon_{k}$):

\vspace*{0.3cm}
{\bf Algorithm~4}
\begin{itemize}
\item{Given: $u_{0}(n,\alpha)$, $t$, and $\epsilon$.}
\item{For $k=0,1,2,\dots$, while $\epsilon_{k} > \epsilon$, do:
\begin{align*}
&u_{k+1}(n,\alpha) = \frac{u_{0}(n(\alpha\delta\, t u_{k}(n,\alpha)-1)^{2},\alpha)}
{\alpha\delta\, t\, u_{0}(n(\alpha\delta\, t u_{k}(n\alpha)-1)^{2},\alpha) + 1},
\\
&\epsilon_{k}=
\left\Vert u_{k+1}(n,\alpha) - u_{k}(n,\alpha)\right\Vert.
\end{align*}
}
\end{itemize}
\begin{figure}[t]
\begin{center}
\hspace*{-0.5cm}
\includegraphics[width=13cm]{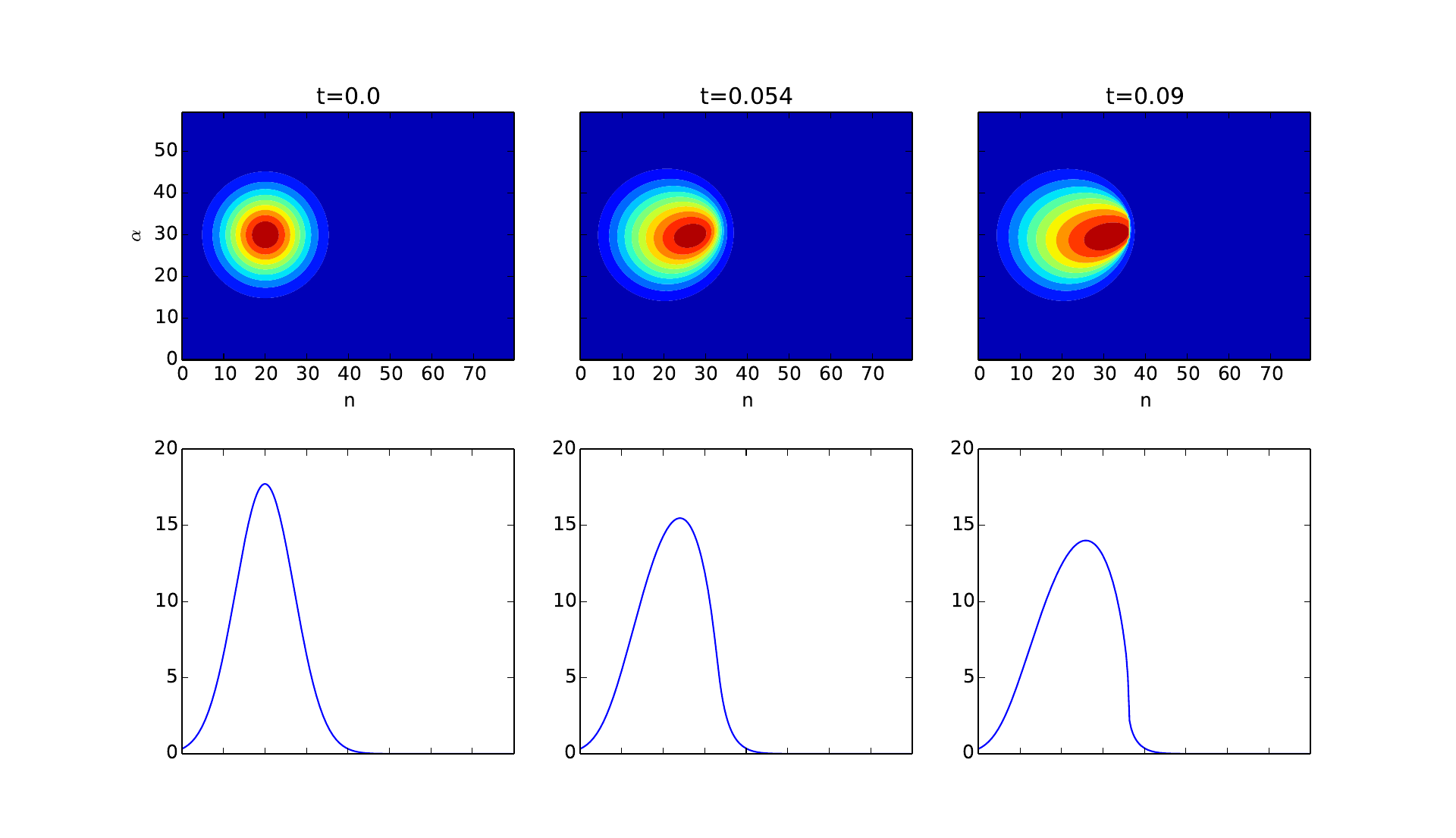}
\end{center}
\caption{
\label{fig:Fig7}
Time evolution of the distribution function $u(n,\alpha,t)$ for populations 
with growth rates depending on $u$, see eq.~(\ref{eq:uNU1}). 
Top row -- snapshots at $t=0.0, 0.054$, and $0.09$; Bottom 
row -- corresponding numerically integrated $\rho(n,t)$ showing the distribution of
populations over their sizes.
}
\end{figure}

Figure~\ref{fig:Fig7} illustrates the application of this algorithm 
(with $\epsilon=10^{-4}$ and $\delta=0.1$).
The initial distribution is the same as in the example of the previous section.
As expected from the non-viscous 
Burgers equation the solution develops a shock (notice the almost vertical front in
the plot at the bottom-right of Fig.~\ref{fig:Fig7}). Even though we have used the 
distribution computed at the previous time instant 
as the initial guess for the next time instant, the convergence of the 
Algorithm~4 slows down as time grows. Around $t=0.1$, as one approachers the 
shock in this example, the convergence of the Algorithm~4 breaks down completely: 
initially the error $\epsilon_{k}$ decreases, however, it never reaches the tolerance 
of $\epsilon=10^{-4}$ and even begins to increase as iterations continue.

Obviously, the development of shock and the related breakdown of the Algorithm~4
does not mean that the populations stop growing. Rather, it is a result of the midpoint 
approximation of the integral in the equation (\ref{eq:DynamicNU}) that reduces the 
problem to the Burgers equation with its well-known shock behavior. Keeping that in mind, 
we see that at least in the early stages of evolution there will emerge a sharp upper 
bound on the sizes of populations and many populations with the average growth factor 
will have that size.

\section{Conclusions and possible applications}
The phase-space approach presented above 
allows computing the evolution of the distribution function of simultaneously 
growing and interacting populations by solving a low-dimensional PDE. 
Although, such phase spaces
impose substantial constraints on the types of interactions, many practical 
problems appear to satisfy these constraints.
As the examples demonstrate, the main power of this approach is that, unlike 
typical Monte-Carlo simulations, it often provides rather
explicit results that could be used as benchmark solutions.

In the case of precision agriculture (Hautala and Hakoj{\"a}rvi~2011), 
this method could be employed to predict and, perhaps, control the distribution 
of plant sizes at the time of harvest. To this end, models considered above 
can be easily extended to include deterministic (measured) 
time-varying growth coefficients and various environmental factors.
Another straightforward extension could be the inclusion of the source term in 
the conservation law (\ref{eq:GenContinuity}) reflecting the emergence of new 
populations, e.g., potato tubers. 

The case of populations competing for limited resources could serve as a 
mathematical model for recovering the distribution of populations from the resource
consumption curves. For example, in the analysis of seed quality it is common to 
monitor the oxygen consumption during the germination process (Bewley et~al.~2013, 
van~Duijn and Koenig~2001, van~Asbrouck and Taridno~2009, Budko et al. 2013), 
with the main consumers of oxygen being the growing populations of 
mitochondria present in every seed.
Our method provides the necessary link between the distribution of the seed 
parameters, including the initial number of active mitochondria,
and the measured oxygen uptake curves for a batch of seeds germinating in a 
closed container. 

On the theoretical level, our phase-space analysis, although different in form from 
the traditional approaches (Barab{\'a}s and Mesz{\'e}na~2009,
Gyllenberg and Mesz{\'e}na 2005), confirms the non-existence of a stable equilibrium for 
coexisting species featuring a continuum of survival strategies, if these strategies
involve the range of carrying capacities. On the other hand, we show that
an asymptotically stable equilibrium is possible, if the inter-species competition 
concerns other survival strategies, such as the rates of growth or 
resource consumption.

Low-dimensional conservation laws occur in migration studies as well. In particular,
we have demonstrated that a completely randomized migration, corresponding to the case
of spatial super diffusion, leads to an explicit phase-space solution resembling 
a time-dependent mollifier of the Dirac's delta function. It is also possible to
incorporate simple models of migration bias in this formulation, as long as these
models do not involve any specific spatial ordering.

Finally, the phase-space coupling model considered in the last section
may have potential applications in the analysis of dynamic systems, where the control
is driven by the time-varying statistical data. For instance, in economics 
and social sciences, the decisions are often based on the perceived actions 
of the majority or minority.

It is easy to imagine that the vector phase-space dynamics
featuring phase-space currents with several nonzero components is much more exciting.
Such dynamics is expected if the growth coefficients or environmental factors 
are varying in time or with the classical predator-prey problem. 
The corresponding continuity equations, however, would
have to be solved numerically (Leveque~2002). 

\section*{Acknowledgements}
The authors are grateful to Prof. A.~van~Duijn (Leiden University and Fytagoras BV),
Dr. S.~Hille (Leiden University),  
Dr.~H.~van Doorn (HZPC Holland BV), and the organizers 
of the European Study Group Mathematics With Industry 
(2013 in Leiden and 2014 in Delft) for the introduction into 
the problems of seed germination and precision agriculture. 


\end{document}